\documentclass[12pt]{article}

\usepackage{amsfonts}
\usepackage[centertags]{amsmath}
\usepackage{amssymb}
\usepackage{amsthm}
\usepackage{newlfont}
\usepackage{graphics}
\usepackage{graphicx}
\usepackage[all]{xy}
\usepackage[german]{varioref}
\usepackage{a4,latexsym}
\newtheorem{thm}{Theorem}
\newtheorem{prop}{Proposition}[section]
\newtheorem{cor}[prop]{Corollary}

\newtheorem{defn}[prop]{Definition}
\newtheorem{rem}[prop]{Remark}
\newtheorem{ex}[prop]{Example}

\newcounter{diagramm}

\newcommand{\fs}{\footnotesize}
\newcommand{\scs}{\scriptsize}

\newcommand{\CC}{\mathbb{C}}

\newcommand{\RR}{\mathbb{R}}
\newcommand{\ZZ}{\mathbb{Z}}
\newcommand{\NN}{\mathbb{N}}
\newcommand{\HH}{\mathbb{H}}

\newcommand{\GGG}{\mathcal{G}}

\newcommand{\aut}{\mbox{Aut}}
\newcommand{\autplus}{\aut^+}
\newcommand{\emaut}{\mbox{\em Aut}}  
\newcommand{\emautplus}{\mbox{\em Aut}^+}

\newcommand{\id}{\mbox{id}}

\newcommand{\affplus}{\mbox{Aff}^+}
\newcommand{\emaffplus}{\mbox{\em Aff}^+}

\newcommand{\veech}{\Gamma}

\newcommand{\slzwei}{\mbox{SL}_2}
\newcommand{\glzwei}{\mbox{GL}_2}

\newcommand{\glk}{\mbox{GL}_k} % fuer lin.tex
\newcommand{\emglk}{\mbox{\em GL}_k}
\newcommand{\emslzwei}{\mbox{\em SL}_2}

\newcommand{\out}{\mbox{Out}}

\newcommand{\stab}{\mbox{Stab}}
\newcommand{\stabaut}{\mbox{Stab}_{\mbox{\fs Aut}^+(\fs F_2)}}

\newcommand{\emstab}{\mbox{\em Stab}}
\newcommand{\emstabaut}{\mbox{\em Stab}_{\mbox{\em \small Aut}^+(\fs F_2)}}

\newcommand{\Xstern}{X^*}

\newcommand{\Estern}{E^*}

\newcommand{\zn}{\ZZ/n\ZZ}

\newcommand{\ifff}{\Leftrightarrow}

\newcommand{\Hquer}{\overline{H}}

\newcommand{\betahat}{\hat{\beta}}

\newcommand{\norm}{\mbox{Norm}}

\newcommand{\mg}{M_g}
\newcommand{\tg}{T_g}

\newcommand{\bpm}{\begin{pmatrix}}
\newcommand{\epm}{\end{pmatrix}}
\newcommand{\re}{\mbox{Re}}
\newcommand{\im}{\mbox{Im}}

\begin{document}

\begin{center}
\section*{Origamis with non congruence Veech groups}
{\it Gabriela Schmith\"usen}
\end{center}

\thispagestyle{plain}
In this article we give an introduction to origamis (often also called 
square-tiled surfaces) and their Veech groups. As main theorem we prove 
that  in each genus there exist origamis, whose Veech groups 
are non congruence subgroups of $\slzwei(\ZZ)$.\\

The basic idea of an {\em origami} is to obtain a 
topological surface  from a few combinatorial data
by gluing finitely many Euclidean unit squares
according to specified rules. 
These surfaces come with a natural 
{\em translation structure}. 
One assigns in general to a translation surface a subgroup of 
GL$_2(\RR)$ called {\em the Veech group}. In the case of surfaces defined 
by origamis, the
Veech groups are finite index subgroups of $\slzwei(\ZZ)$. These groups are
the objects we study in this article.\\

One motivation to be interested in Veech groups is their relation 
to 
{\em Teichm\"uller disks} and {\em Teichm\"uller curves}, see e.g. the article
\cite{h-yokohama} of F. Herrlich in the same volume:
A translation surface of genus $g$ defines in a geometric way an embedding of 
the upper half plane into the Teichm\"uller space $T_g$ of 
closed Riemann surfaces of genus $g$. The image is called Teichm\"uller disk.
Its projection to the moduli space $M_g$ 
is sometimes a complex algebraic curve, called
{\em Teichm\"uller curve}. 
More precisely this happens, if and only if the Veech group is a lattice in 
$\slzwei(\RR)$. In this case the algebraic curve can be determined from
the Veech group up to birationality.\\

It is hard to determine the Veech group for a general translation surface.
However, if the translation surface comes from an origami there is a special
approach to this problem. It is based on the idea of describing origamis 
by finite index subgroups of $F_2$, the free group in two generators. 
This leads to a characterization
of origami Veech groups as the images in $\slzwei(\ZZ)$ of certain subgroups of $\aut(F_2)$, the automorphism group of $F_2$.\\

Using this approach we will calculate Veech groups of two origamis 
explicitly. They turn out to be non congruence groups. Starting from these
examples we obtain infinite sequences of origamis all of whose Veech groups are
non congruence groups. This leads to the following theorem.
\begin{thm} \label{main}
Each moduli space $M_g$  $(g \geq 2)$ contains an origami curve whose Veech
group is a non congruence group.
\end{thm}

In Section \ref{origamis} we introduce origamis and present
different equivalent ways to describe them. In Section 
\ref{veechgroups}
we give a glance on the mathematical context. 
We describe, how an origami defines a family of 
translation surfaces and  
explain roughly , how one obtains  a Teichm\"uller curve in moduli space
starting from an 
origami. We introduce Veech groups
and shortly point out their relation to {Teich}\-m\"uller curves. In 
Section
\ref{ovgs} we turn to Veech groups of origamis and present
a characterization of them in terms of 
automorphisms of the free group $F_2$ in two generators. We use this
characterization to calculate two examples explicitly. Finally,
in Section \ref{noncgsection} we show that these two examples
produce Veech groups that are non congruence groups and give
a method to construct out of them infinite sequences of Veech groups that are
again non congruence groups.\\

The first part (Section \ref{origamis} - Section \ref{ovgs}) 
of this article is meant to give a handy introduction 
to origamis and an overview on some of our results about their Veech groups. 
In the second part we state and prove Theorem \ref{main} based
on the results in the PhD thesis \cite{gabidiss} of the author.\\

For a broader introduction and overview on
origamis and Teichm\"uller curves as well as for references to the
larger context, 
we refer the the reader e.g. to \cite{hesh2}, \cite{sh1} and \cite{gabidiss}.\\

{\bf Acknowledgments:} I would like to thank Frank Herrlich
for his support in respect of the content and for his proof reading,
Stefan K\"uhnlein for  helpful discussions and suggestions especially on
non congruence groups and
the organizers of the conference for giving me the opportunity to
contribute to these proceedings.
This work was partially supported by a fellowship within the Postdoc-Programme
of the German Academic Exchange Service (DAAD).

\newpage
\section{Origamis}\label{origamis}

There are several ways to define origamis. We start with the 
somehow playful description that we have learned from \cite{l}, where also 
the name
{\em origami} was introduced:
An {\em origami} is obtained by 
gluing the edges of finitely many copies $Q_1$, \ldots, $Q_d$ 
of the Euclidean square $Q$ via translations according to the
following rules:
\begin{itemize}
\item
Each left edge shall be identified to a right edge and vice versa.
\item
Similarly, each upper edge shall be identified 
to a lower one.
\item The arising closed surface $X$ shall be  connected.
\end{itemize}

We only study what is called oriented origamis in \cite{l} and call them
just {\em origamis}.

\begin{ex}\label{ex1}\hspace{1cm}\\[-7mm]
\begin{enumerate}
\item[a)] 
The simplest example is the origami that is made from only one square.
There is precisely one possibility to glue its edges according to the 
rules.
One obtains a torus $E$. We call this origami the {\em trivial origami}
$O_0$.

\setlength{\unitlength}{1cm}
\begin{center}
\begin{picture}(1.1,1.1)
\put(0,0){\framebox(1,1){}}
% Labels:
\put(.4,-.35){$a$}
\put(1.1,.4){$b$}
\put(.4,1.1){$a$}
\put(-.3,.4){$b$}
\put(-.4,-.4){$\infty$}
\put(-.13,-.13){\Large $\bullet$}
\put(-.13, .87){\Large $\bullet$}
\put(.87,.87){\Large $\bullet$}
\put( .87,-.13){\Large $\bullet$}
\end{picture}\\[3mm]
\end{center}
\begin{center}
\refstepcounter{diagramm}{\it Figure \arabic{diagramm}}:
{\it The trivial origami. Opposite edges are glued.}
\label{trivial}
\end{center}

Observe that the four vertices of the square are all identified
and become one point on the closed surface $E$. We call 
this point $\infty$.

\item[b)]
We now consider an origami made from four squares, see Figure~\ref{l23}.
Some identifications of the edges are already done in the picture. For
all other edges those having same labels are glued. The origami
is called  $L(2,3)$ for obvious reasons.

\setlength{\unitlength}{1cm}
\begin{center}
\begin{picture}(3,2.3)
\put(0,0){\framebox(1,1){2}}
\put(1,0){\framebox(1,1){3}}
\put(2,0){\framebox(1,1){4}}
\put(0,1){\framebox(1,1){1}}

% Beschriftung:
\put(-.13,-.13){\Large $\bullet$}
\put(-.13, .87){\Large $\bullet$}
\put(-.13,1.87){\Large $\bullet$}
\put( .87,-.13){\Large $\bullet$}
\put( .87, .87){\Large $\bullet$}
\put( .87,1.87){\Large $\bullet$}
\put(2.87,-.13){\Large $\bullet$}
\put(2.87, .87){\Large $\bullet$}

\put(1.87,-.14){\Large\bf $\circ$}
\put(1.87, .84){\Large $\circ$}

%\put(1.75,-.5){\LARGE \bf *}
%\put(1.75, .45){\LARGE \bf *}

% Labels:
\put(.4,-.35){$a$}
\put(1.5,-.35){$b$}
\put(2.5,-.35){$c$}
\put(3.1,.4){$d$}
\put(2.5,1.1){$c$}
\put(1.5,1.1){$b$}
\put(.4,2.1){$a$}
\put(-.3,.4){$d$}
\put(1.1,1.45){$e$}
\put(-.3,1.45){$e$}

\end{picture}\\[3mm]
\end{center}
\begin{center}
\refstepcounter{diagramm}{\it Figure \arabic{diagramm}}:
{\it The origami $L(2,3)$. Opposite edges are glued.}
\label{l23}
\end{center}

Observe that in this case the vertices labeled with $\bullet$
and the vertices labeled with $\circ$ are
respectively identified and become two
points  on the closed surface $X$. 
By calculating the  Euler characteristic
one obtains, that the genus of the surface $X$  is 2.
\item[c)]
Finally, we consider an example with five squares, see Figure \ref{d}.
Here, edges with same labels are identified. For the unlabeled edges, those
which are opposite to each other are glued. We call the origami $D$.
\begin{center}
\begin{picture}(3,3.3)
\setlength{\unitlength}{1cm}
\begin{picture}(4,3.2)
% Die Kaestchen des Origamis:
%\put(-2.7,1.2){ $\mbox{Loc}(n) =$}
\put(0,0){\framebox(1,1){1}}
\put(1,0){\framebox(1,1){2}}
\put(2,0){\framebox(1,1){3}}
\put(0,1){\framebox(1,1){4}}
\put(0,2){\framebox(1,1){5}}

% Beschriftung:
\put(1.55,-.35){$a$}
\put(2.55,1.1){$a$}
\put(1.55,1.1){$b$}
\put(2.55,-.35){$b$}

% labels for the vertices:
\put(-.1,-.1){\Large $\bullet$}
\put( .9,-.1){\Large $\bullet$}
\put(2.9,-.1){\Large $\bullet$}
\put(1.9, .9){\Large $\bullet$}
\put(-.1,2.9){\Large $\bullet$}
\put( .9,2.9){\Large $\bullet$}
\put(-.13, 1.87){\Large $\circ$}
\put( .87, 1.87){\Large $\circ$}
\put(-.25,  .7){\LARGE \bf *}
\put( .75,  .7){\LARGE \bf *}
\put(1.7,  -.3){\LARGE \bf *}
\put(2.7,  .7){\LARGE \bf *}
\end{picture}

\end{picture}\\[3mm]
\end{center}
\begin{center}
\refstepcounter{diagramm}{\it Figure \arabic{diagramm}}:
{\it The origami $D$. Edges with the same label and unlabeled
edges that are opposite are glued.}
\label{d}
\end{center}
In this case, we obtain  
the three identification classes
$\circ$, $\star$ and $\bullet$ for the vertices. 
The genus of the closed surface $X$ is
again 2.
\end{enumerate}
\end{ex}

\subsubsection*{Origamis as coverings of a torus}

Observe, that the trivial origami $O_0$ from Example \ref{ex1}~a) is 
universal in the following
sense: If $X$ is the closed surface that arises from an arbitrary
origami $O$ and $E$ the torus that arises from $O_0$, then
we have a natural map $X \to E$ by mapping each of the unit squares
of the origami $O$ that form the surface $X$ to the one unit square of 
$O_0$ 
that forms the torus $E$.
This map is a covering that is unramified except over the one point 
$\infty \in E$. Conversely, given a closed surface $X$ together with 
such a covering $p:X \to E$, we obtain a decomposition of $X$ into squares 
by cutting $X$ along the preimages of the edges of the one square of $O_0$ 
that forms $E$.
This motivates the following definition of {\em origamis}.
\begin{defn}
An {\em origami} $O$ of {\em genus} $g$ and {\em degree} $d$ 
is a covering $p:X \to E$ of degree $d$ from a  
closed, oriented (topological) surface $X$ of genus $g$
to the torus $E$ that is ramified over at most one marked point $\infty \in E$.
\end{defn}

Remark that we have fixed here one torus $E$ and one point $\infty \in E$.
In particular we may furthermore fix a point $M \neq \infty$ on $E$
and a set of standard generators 
of the fundamental group $\pi_1(E,M)$ that do not pass through $\infty$. 
That way we obtain a fixed  isomorphism 
\begin{equation}\label{isom}
\pi_1(\Estern) \cong F_2,
\end{equation}
where $\Estern = E - \infty$ and $F_2 = F_2(x,y)$ is the free group in 
two generators $x$ and $y$.
Describing E by gluing the edges of the unit square 
via translations, we choose $M$ to be the midpoint of the unit square 
and the standard generators to be the horizontal and the vertical
simply closed curve through $M$, see Figure \ref{xybild}.

\setlength{\unitlength}{1.5cm}
\begin{center}
\begin{picture}(1,1)
\put(0,0){\framebox(1,1){}}
\put(0,.5){\vector(1,0){1}}
\put(.5,0){\vector(0,1){1}}
\put(.45,.45){$\bullet$}
\put(.55,.55){$M$}
\put(.15,.55){$x$}
\put(.55,.15){$y$}
\end{picture}\\[3mm]
\refstepcounter{diagramm}{\it Figure
\arabic{diagramm}\label{xybild}: Generators of $\pi_1(\Estern)$.}
\end{center}

\begin{ex}
In Example \ref{ex1}, in a) the 
covering is the identity $\id: E \to E$.\\ 
In b) we have a covering 
$p:X \to E$ of degree 4 that is ramified in the two points 
labeled by $\bullet$ and $\circ$. Recall that the genus of $X$ is 2.\\ 
In c) we have
a covering $p:X \to E$ of degree 5 ramified in the two points 
labeled by $\bullet$ and
$\star$. Observe that though the point on $X$ labeled by $\circ$ is 
a preimage of $\infty$, the covering is not ramified in this point. The genus
of $X$ is again 2.
\end{ex}

\begin{defn}
We say that two origamis $O_1 = (p_1:X_1 \to E)$ and
$O_2 = (p_2:X_2 \to E)$ are {\em equivalent}, if there is a homeomorphism
$\varphi: X_1 \to X_2$ with $ p_1 = p_2 \circ \varphi$.
\end{defn}

\subsubsection*{Description by a pair of permutations}

An origami $O = p:(X \to E)$ of degree $d$ defines 
(up to conjugation in $S_d$)
\begin{itemize}
\item
a homomorphism
$m:F_2=F_2(x,y) \to S_d$  or equivalently
\item
a pair of permutations $(\sigma_a,\sigma_b)$ in $S_d$  
\end{itemize}
as follows:\\
Let $M_1$, \ldots, $M_d$ be the preimages of the point $M$ (defined
as above) under $p$. Furthermore, let
\[m:\pi_1(\Estern,M) \to \mbox{Sym}(M_1,\ldots,M_d)\] be the {\em 
monodromy map} 
defined by
$p$, i.e. for the closed path $c \in \pi_1(\Estern,M)$ the point
$M_i$ is mapped to $M_j$ by $m(c)$ if and only if the lift of the
curve $c$ to $X$ via $p$, that starts in $M_i$, ends in $M_j$.\\
Choosing an isomorphism Sym$(M_1,\ldots, M_d) \cong S_d$
and using the isomorphism $\pi_1(\Estern) \cong F_2$ fixed in  (\ref{isom})
makes $m$ into a homomorphism from $F_2$ to $S_d$. We set
$\sigma_a = m(x)$ and $\sigma_b = m(y)$.\\
Observe that this homomorphism depends on the chosen isomorphism to $S_d$ and 
on the
choice of the origami in its equivalence class only up to
conjugation in $S_d$. Therefore we consider two homomorphisms 
$m_1: F_2 \to S_d$ and $m_2: F_2 \to S_d$ to be {\em equivalent}, if 
they are conjugated by an element in $S_d$. Similarly we call
two pairs $(\sigma_a,\sigma_b)$ and $(\sigma_a',\sigma_b')$ in $S_d$
{\em equivalent}, if they are simultaneously conjugated, i.e.
there is some $s \in S_d$ such that $\sigma_a = s\sigma_a's^{-1}$
and $\sigma_b = s\sigma_b's^{-1}$.

\begin{ex}
In Example \ref{ex1} we obtain for the origami $L(2,3)$ in b)
the monodromy homomorphism
\[m: F_2 \to S_4, \quad 
       x  \mapsto (2\ 3\ 4) \; \mbox{ and } \; y \mapsto (2\ 1),\] 
and thus $\sigma_a = (2\ 3\ 4)$ and $\sigma_b = (2\ 1)$.\\
For the origami $D$ in c) we similarly obtain
the permutations 
\[\sigma_a = (1\ 2\ 3) \;\; \mbox{ and } \;\; \sigma_b = (1\ 4 \ 5)(2\ 3).\]

\end{ex}

\subsubsection*{Description as finite index subgroups of \boldmath{$F_2$}}

Origamis can be equivalently described as finite index subgroups of $F_2$, the 
free group in two generators, as stated in the following remark. 
The characterization of the Veech groups of origamis is mainly based on this 
observation.

\begin{rem}\label{corespond}
We have a one-to-one correspondence:
\[
\mbox{origamis up to equivalence} \;\; \leftrightarrow \;\;
\mbox{finite index subgroups
of $F_2$ up to conjugacy}.
\]
\end{rem}

More precisely, this correspondence is given as follows:\\
Let $O = (p:X\to E)$ be an origami. Define $\Estern = E - \{\infty\}$ 
and $\Xstern = X - p^{-1}(\infty)$. Thus we may restrict $p$ to the unramified 
covering $p:X^* \to E^*$. This defines an embedding of the corresponding
fundamental groups:
\[U = \pi_1(\Xstern) \;\; \subseteq \;\; \pi_1(\Estern) \cong F_2\]
Again we use the fixed isomorphism in (\ref{isom}), see also Figure 
\ref{xybild}.
Changing the origami in its equivalence class leads to a conjugation
of $U$ with an element in $F_2$. The index of the subgroup of $F_2$
is the degree $d$ of the covering $p$.\\ 
Conversely, given a finite index subgroup $U$ of $F_2$ we retrieve the
origami in the following way:
Let $v: \tilde{\Estern} \to \Estern$ be a universal covering of $\Estern$.
By the theorem of the universal covering,  
$\pi_1(\Estern)$ is isomorphic to Deck$(\tilde{\Estern}/\Estern)$, 
the group of deck transformations of $\tilde{\Estern}/\Estern$. Furthermore,
the finite index subgroup $U$ of Deck$(\tilde{\Estern}/\Estern)$ corresponds to an 
unramified covering $p:\Xstern \to \Estern$ of finite degree. This can be extended
to a covering $X \to E$, where $X$ is a closed surface.

\begin{ex}\label{fgex}
In Example \ref{ex1}, we obtain the following subgroups of $F_2$:\\
In a), $\Xstern$ is the once punctured torus itself and $U = F_2$.\\
In b), $\Xstern$ is a genus 2 surface with 2 punctures. Thus 
$U = \pi_1(\Xstern)$ is a free group of rank 5. Keeping in mind that 
we use the identification $\pi_1(\Estern) \cong F_2 = F_2(x,y)$ 
described in Figure \ref{xybild}, one can read off from the
picture in Figure \ref{l23} that 
\[U \; = \;\; <x^3,\; xyx^{-1},\; x^2yx^{-2},\; yxy^{-1},\; y^2>\]
In c), $\Xstern$ is a genus 2 surface with three punctures. Thus $U$
is a free group of rank 6. More precisely, we read off the picture in  Figure \ref{d},
that
\[ U \; = \;\;  <x^3,\; xyx^{-2},\; x^2yx^{-1},\; yxy^{-1},\; y^2xy^{-2},\; y^3>\]
\end{ex}

\subsubsection*{Description as a finite graph}

Finitely, sometimes it is convenient to describe 
an origami $O = (p:X \to E)$ as a finite, oriented
labeled graph: Namely, let $U$ be the finite index subgroup of $F_2$
(unique up to conjugation) that corresponds to $O$ as described in the
last paragraph. Then we represent the origami by the Cayley-Graph of 
$U\subseteq F_2$: The vertices of the graph are the coset representatives.
They are labeled with a representative of the coset. The edges are labeled
with $x$ and $y$. For each
vertex (with label $w \in F_2$) there  is an $x$-edge from it
to the vertex that belongs to the coset of $wx$. And similarly there is
a $y$-edge to the vertex that belongs to the coset $wy$.

\begin{ex}
The following figure shows the Cayley-graph for the origami $L(2,3)$ 
from Example \ref{ex1}:\\

\begin{minipage}{\linewidth}
\[
%\SelectTips{cm}{}
\xymatrix @-1pc {
*++[o][F-]{\bar{y}}\ar@/_/[d]_y \ar@(dr,ur)[]_{x}&&\\
*++[o][F-]{\bar{\mbox{id}}}\ar[u]_y\ar[r]^{x}&
*++[o][F-]{\bar{x}}\ar[r]^{x}\ar@(ur,ul)[]_y&
*++[o][F-]{\bar{x^2}}\ar@(ur,ul)[]_y \ar@/^1pc/[ll]^x
}
\]
\makebox{
\begin{minipage}{\linewidth}
\begin{center}
\refstepcounter{diagramm}
{\it Figure }\arabic{diagramm}: Graph for $O = L(2,3)$.\label{l23tree}
\end{center}
\end{minipage}
}
\end{minipage}
\end{ex}

\section{Translation structures and Veech groups}\label{veechgroups}

\subsubsection*{Translation structures}

Recall that an atlas on a surface is called 
{\em translation atlas}, if all transition maps are translations.
An origami  $O=(p:X \to E)$ naturally defines
an SL$_2(\RR)$-family of translation
structures $\mu_A$ ($A \in \slzwei(\RR)$) 
on $\Xstern = X - p^{-1}(\infty)$ as follows:
\begin{itemize}
\item
As first step, observe that each $A \in \slzwei(\RR)$
naturally defines a translation structure $\eta_A$ on the torus $E$ itself 
by identifying it with $\CC/\Lambda_A$, where
\begin{equation} \label{lattice}
A = \bpm a&b\\c&d \epm \;\; \mbox{ and } \;\;
\Lambda_A \mbox{ is the lattice } <\bpm a\\c\epm,\bpm b\\d \epm>
\mbox{ in } \CC 
\end{equation}
\item Then define the translation structure $\mu_A$ on $\Xstern$
by lifting $\eta_A$ via $p$, i.e. $$\mu_A = p^*\eta_A.$$
\end{itemize}
Using the first description of an origami that we gave  by gluing 
squares, we obtain the translation structure $\mu_I$ (where $I$
is the identity matrix), if we identify the squares with the Euclidean
unit square in $\CC$. We obtain $\mu_A$ for a general matrix 
$A \in \slzwei(\RR)$ from this by identifying the squares with the 
parallelogram spanned by the two vectors
\[\bpm a\\c \epm \mbox{ and } \bpm b\\d \epm .\]
Thus the $\slzwei(\RR)$-variations of the translation
structure $\mu_I$ can be thought of as affine shearing of the unit
squares, see Figure \ref{schraeg}.\\[2mm]
\setlength{\unitlength}{1cm}
\hspace*{3.5cm}
\begin{minipage}{6cm}
\begin{center}
\begin{picture}(4,2.5)
\put(0,0){\line(1,1){1}}
\put(0,0){\line(1,0){1}}
\put(1,1){\line(1,0){1}}
\put(1,0){\line(1,1){1}}
\put(1,0,){\line(1,0){1}}
\put(2,0,){\line(1,0){1}}
\put(2,1,){\line(1,0){1}}
\put(2,2){\line(1,0){1}}
\put(2,2){\line(1,0){1}}
\put(2,0,){\line(1,1){1}}
\put(2,1,){\line(1,1){1}}
\put(2,1,){\line(1,1){1}}
\put(1,1){\line(1,1){1}}
\put(3,1){\line(1,0){1}}
\put(3,0){\line(1,1){1}}
\end{picture}

\end{center}
\end{minipage}\\[3mm]
\begin{center}
\refstepcounter{diagramm}{\it Figure \arabic{diagramm}}:
{\it Sheared translation structure for the origami $L(2,3)$.}
\label{schraeg}
\end{center}

\subsubsection*{From an origami to a Teichm\"uller curve in the moduli space}

By the $\slzwei(\RR)$-family of translation structures, the origami
$O = (p:X \to E)$ defines a specific 
complex algebraic curve called {\em Teichm\"uller curve}
in the moduli space $\mg$ of closed Riemann
surfaces of genus $g$. We state this construction here only briefly as 
motivation and refer e.g. to the overview article \cite{hesh2} for a 
detailed description and links to references. A particular nice
configuration of such Teichm\"uller curves is described in 
\cite{h-yokohama} in this volume.\\ 

The Teichm\"uller curve in $\mg$ 
is obtained from the origami in the following way:
\begin{itemize}
\item The translation structure $\mu_A$ described in the previous paragraph
is in particular a complex structure on the surface
$\Xstern$ which can be extended to the closed surface $X$. The Riemann surface
$(X,\mu_A)$ together with the identity map $\id:X \to X$ as marking then 
defines
a point in the Teichm\"uller space $\tg$. Thus we obtain the map: \; 
$\hat{\iota}:\;\slzwei(\RR) \to \tg, \; A \mapsto [(X,\mu_A),\id]$.
\item
If $A \in \mbox{SO}_2(\RR)$, then the affine map $z \mapsto A\cdot z$
is holomorphic. Thus the map $\hat{\iota}$ factors through 
SO$_2(\RR)$. Furthermore
using that $\slzwei(\RR)$ modulo SO$_2(\RR)$  is isomorphic to the
upper half plane $\HH$, one obtains a map
\[\iota:\; \HH \cong \mbox{SO}_2(\RR)\backslash\slzwei(\RR) 
\;\;\; \to \;\; \tg\]
In fact, this map is an embedding that is in the same time holomorphic and
isometric. A map with this property  is called {\em Teichm\"uller embedding}
and its image $\Delta$ 
in Teichm\"uller space is called a {\em Teichm\"uller disk} or
a {\em geodesic disk}.
\item
Finally, one may compose the map $\iota$ with  the projection to the moduli space
$\mg$. The image of $\Delta$ in $\mg$ is a complex algebraic curve. 
A curve in $\mg$ that arises like this as the image of a 
Teichm\"uller disk is called 
{\em Teichm\"uller curve}.
\end{itemize}

{\bf Note:} More generally, one obtains a Teichm\"uller disk $\Delta$
in a similar way starting from an arbitrary translation surface (or a bit more
general: from a flat surface). However, the image of such a disk $\Delta$
in moduli space is not always a complex algebraic curve; in fact its 
Zariski closure tends to be of higher dimension. It is an interesting
question how to decide whether a translation surface leads to a Teichm\"uller 
curve. One possible answer to this is given by the Veech group which we
introduce in the following paragraph.

\subsubsection*{Veech groups}

Let $\Xstern$ be a connected surface and $\mu$ a translation structure
on it. One assigns to it a subgroup of GL$_2(\RR)$ called {\em Veech group}
as described in the following: We consider the group
$\affplus(\Xstern,\mu)$ of all orientation preserving 
{\em affine diffeomorphisms}, i.e. 
orientation preserving diffeomorphisms  
that are locally affine maps of the plane $\CC$, see Figure \ref{affine}. 
Here
-- and throughout the whole article -- we identify $\CC$ with $\RR^2$ by the map 
$z \mapsto (\re(z) \,,\, \im(z))^t$. Thus an affine diffeomorphism $f$
can be written in terms of local charts as
\begin{equation}
f:\; z = (\re(z),\im(z))^t \mapsto A\cdot (\re(z),\im(z))^t + z_0
\; \mbox{ with } A \in \glzwei(\RR) \mbox{ and } z_0 \in \CC.
\end{equation}

Observe that $A$ does not depend on the chart, since $\mu$ is
a translation structure. Thus one obtains a well defined map 
\[D:\, \affplus(\Xstern,\mu) \to \glzwei(\RR),\;\; f \mapsto A\]
called {\em Derivative map}.
\begin{defn}
The {\em Veech group} $\Gamma(\Xstern,\mu)$ of the translation 
surface $(\Xstern,\mu)$ is the image of the derivative map $D$:
\[\Gamma(\Xstern,\mu) = D(\emaffplus(\Xstern,\mu))\]
\end{defn}

\vspace*{1mm}

\hspace*{5mm}
\begin{minipage}{5cm}
  \includegraphics[scale=.7]{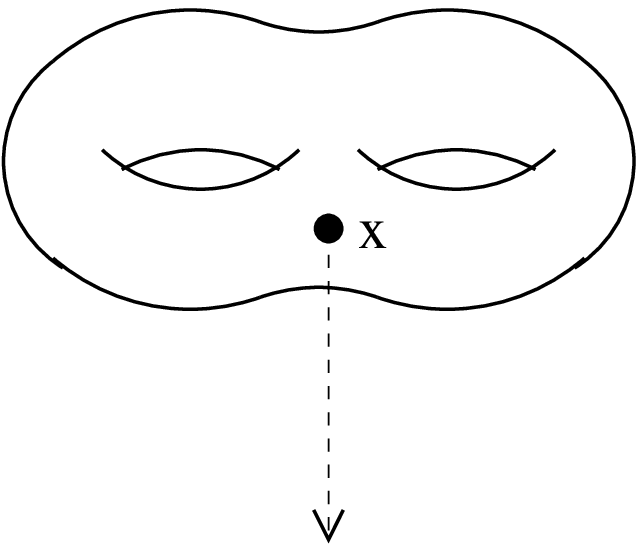}
\end{minipage}
\begin{minipage}{3cm}
 \vspace*{-30mm}
 \hspace*{10mm}   f\\
 \hspace*{-2mm}
  \includegraphics{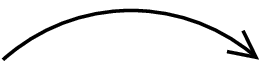}
\end{minipage}
\hspace*{2mm}
  %Flaeche dazu:
\begin{minipage}{7cm}
  \includegraphics[scale=.7]{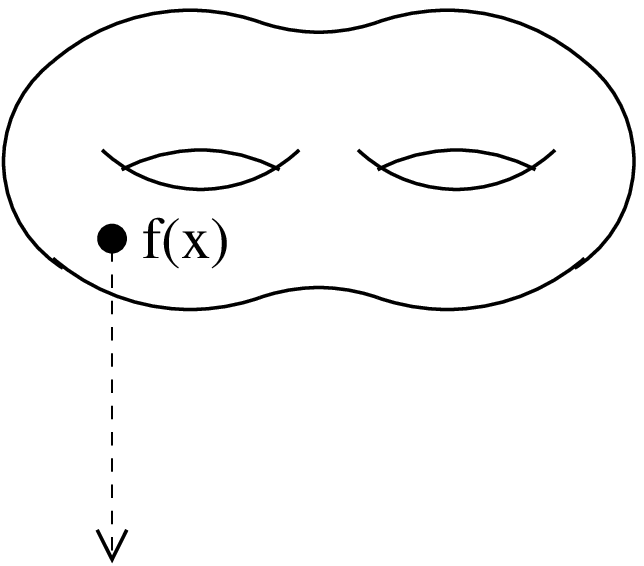}
\end{minipage}\\[-7mm]
\hspace*{2.6cm}
$\CC$ \hspace*{2.3cm}
\begin{minipage}[b]{3cm}
\hspace*{3mm}$z \mapsto Az + z_0$\\[1mm]
\includegraphics{pfeil.eps} 
\end{minipage}
\hspace*{10mm}$\CC$\\[1mm]
\hspace*{15mm}
\refstepcounter{diagramm}{\it Figure \arabic{diagramm}}:
{\it An affine diffeomorphism of a translation surface}
\label{affine}

\begin{ex}\label{torusvg}
Let $(\Xstern,\mu)$ be $\CC/\Lambda_I$ with the natural
translation structure. Here $I$ is the identity matrix
and $\Lambda_I$ is the corresponding lattice as 
defined  in (\ref{lattice}).\\ 
An affine diffeomorphisms of $\CC/\Lambda_I$ lifts to an
affine diffeomorphism of $\CC$ respecting the lattice.
Conversely, each such diffeomorphism descends to $\CC/\Lambda_I$.
Thus, we have in this case
\[\Gamma(\Xstern,\mu) = \emslzwei(\ZZ).\]
\end{ex}

\subsubsection*{Veech groups and Teichm\"uller curves}

As indicated in the paragraph about Teichm\"uller curves,
the Veech group \  ``knows'' whether
a translation surface defines a Teichm\"uller curve in moduli space or not.
More precisely, one has the following statement:\\[2mm]
{\bf Fact:}\label{fact}
Let $X$ be a surface of genus $g$ and $\Xstern = X- \{P_1, \ldots, P_n\}$ 
for finitely many points $P_1$, \ldots, $P_n$ on $X$. Furthermore let $\mu$
be a translation structure on $\Xstern$.\\ 
Then $(\Xstern,\mu)$ defines a
Teichm\"uller curve $C$ if and only if the Veech group $\Gamma(\Xstern,\mu)$
is a lattice in $\slzwei(\RR)$.
In this case, the curve $C$ is (antiholomorphic) 
birational to $\HH/\Gamma(\Xstern,\mu)$.\\

We describe the relation to Teichm\"uller curves 
here just as motivation and in order to give 
a glance at the general frame. We have therefore resumed 
theorems contributed by several authors condensed in what is 
here called ``fact''. 
A good access to it can be found e.g. in \cite{eg} or \cite{z}.
A broader overview on Veech groups of translation surfaces
is given e.g. in \cite{hs} and in \cite{li2}. 
Teichm\"uller disks, Teichm\"uller curves and Veech groups have 
intensively been studied
by numerous authors, 
starting from Thurston \cite{th2} and  Veech himself \cite{ve}. 
We refer to \cite{sh1} and \cite{hesh2} for more comprehensive overviews 
on references.

\section{Veech groups of origamis}\label{ovgs}

Let $O = p:(X \to E)$ be an origami. We have seen in Section \ref{veechgroups}
that $O$ defines an $\slzwei(\RR)$-family of translation structures $\mu_A$
($A \in \slzwei(\RR)$) on $\Xstern = X - p^{-1}(\infty)$. The corresponding Veech groups are not very different. In fact, they are all
conjugated to each other. More precisely, we have:
\[\Gamma(\Xstern,\mu_A) = A\Gamma(\Xstern,\mu_I)A^{-1}.\]
Thus, we may restrict to the case where $A = I$ which justifies the
following definition.

\begin{defn}\label{defvgorigami}
The {\em Veech group $\Gamma(O)$ of the origami $O$} is defined to be
$\Gamma(\Xstern,\mu_I)$. 
\end{defn}

From Example \ref{torusvg} it follows that the Veech group of the trivial
origami $O_0$ (defined in Example \ref{ex1}) is $\slzwei(\ZZ)$. For a 
general origami one can show that $\Gamma(O)$ is a finite index subgroup
of $\slzwei(\ZZ)$. In fact, also the converse is true as it was shown by 
Gutkin and Judge in \cite{gj}:  A Veech group is a finite index subgroup of
$\slzwei(\ZZ)$ if and only if it comes from an origami.\\
From this it follows in particular by the Fact presented in 
Section \ref{veechgroups}
on page \pageref{fact} that an origami always defines a Teichm\"uller curve
in the moduli space.

\subsubsection*{Characterization of origami Veech groups}

Recall from Section \ref{origamis} that an origami $O$ 
corresponds (up to equivalence)
to a finite index subgroup $U$ of $F_2 = F_2(x,y)$, the free
group in two generators (up to conjugation). This description
enables us to give a characterization of its Veech group
entirely in terms of $F_2$ and its automorphisms.\\

For this we need the following two ingredients:
\begin{itemize}
\item Let $\betahat: \aut(F_2) \to \out(F_2) \cong \glzwei(\ZZ)$ 
be the natural projection. The fact that we only consider
orientation preserving diffeomorphisms applies to only taking automorphisms
of $\aut(F_2)$ that are mapped to elements in $\slzwei(\ZZ)$. We denote
$\autplus(F_2) = \betahat^{-1}(\slzwei(\ZZ))$ and restrict to the map
\[\betahat: \autplus(F_2) \to \slzwei(\ZZ).\]
 \item
Let $\stab(U) = \{\gamma \in \autplus(F_2)| \gamma(U) = U \}$
\end{itemize}
Using these ingredients, it was shown in \cite{sh1}
that Veech groups of origamis 
can be described as stated in the following theorem.

\begin{thm}[Proposition 1 in \cite{sh1}] \label{char}
For the Veech group $\Gamma(O)$ of the origami $O$ holds:
\[\Gamma(O) = \betahat(\emstab(U))\]
\end{thm}

Let us make two comments on this description:\\
One consequence is, that 
one obtains an algorithm that can calculate 
the Veech group of an arbitrary origami explicitly. 
This algorithm is described in detail in \cite{sh1}.\\
As an other consequence, we have now a characterization of
all origami Veech groups as stated in the following corollary. 

\begin{cor}
A finite index subgroup of $\emslzwei(\ZZ)$ occurs as origami Veech group 
if and only if 
it is the image of the stabilizing group $\emstab(U) \subseteq 
\emautplus(F_2)$ for some
finite index subgroup $U$ in $F_2$.
\end{cor}

Thus the question, which finite index subgroups of 
$\slzwei(\ZZ)$ are Veech groups
becomes roughly speaking the same as the question which subgroups
of $\autplus(F_2)$ are such stabilizing groups. So far, there is no general
answer known.\\

In \cite{gabidiss} it was shown that many congruence
subgroups of $\slzwei(\ZZ)$ are Veech groups. Recall
that a {\em congruence group} of level $n$ 
is a subgroup of $\slzwei(\ZZ)$
that is the full preimage of some subgroup of $\slzwei(\ZZ/n\ZZ)$
under the natural homomorphism $\slzwei(\ZZ) \to \slzwei(\ZZ/n\ZZ)$
and $n$ shall be minimal with this property.
For prime level congruence groups the following statement is shown
in \cite[Theorem~4]{gabidiss}

\begin{thm}\label{cgs}
Let $p$ be prime. All congruence groups $\Gamma$ of level $p$ are Veech groups
except possibly p $\in \{2,3,5,7,11\}$ and $\Gamma$ has index $p$ in
$\emslzwei(\ZZ)$.
\end{thm}

This result is generalized to a statement for arbitrary $n$ in 
\cite[Theorem 5]{gabidiss}

\subsubsection*{Presenting the Veech group \boldmath{$\Gamma$} 
and the quotient \boldmath{$\HH/\Gamma$}
for an origami}

As mentioned above, using Theorem \ref{char} the Veech group of 
an origami can be  calculated explicitly. The Veech groups
are described as subgroups of $\slzwei(\ZZ)$ by generators and
coset representatives. We use for the notation that 
$\slzwei(\ZZ)$ is generated by $S$ and $T$, with
\[S = \bpm 0&-1\\1&0\epm \;\;\mbox{ and } \;\; 
  T = \bpm1&1\\0&1 \epm.\]
Recall furthermore 
from the discussion on Veech groups and Teichm\"uller curves 
in Section \ref{veechgroups} on page \pageref{fact}
that for a Veech group $\Gamma$ 
we are in particular interested in the quotient $\HH/\Gamma$, 
since this quotient is birational to the corresponding Teichm\"uller curve.
Here $\Gamma$ acts as Fuchsian group on the upper half plane $\HH$, which
is endowed with the Poincar\'e metric.\\
Since an origami Veech group $\Gamma$  is a finite index
subgroup of $\slzwei(\ZZ)$, 
the quotient $\HH/\Gamma$ comes with a natural triangulation. More precisely, 
we choose the fundamental domain for the action of $\slzwei(\ZZ)$ on $\HH$ 
that is the geodesic pseudo-triangle $\Delta$ with vertices 
$P = -\frac{1}{2} + \frac{\sqrt{3}}{2}i$, 
$Q =  \frac{1}{2} + \frac{\sqrt{3}}{2}i$ and 
$P_{\infty} = \infty$.
\begin{center}
\includegraphics[scale = .5]{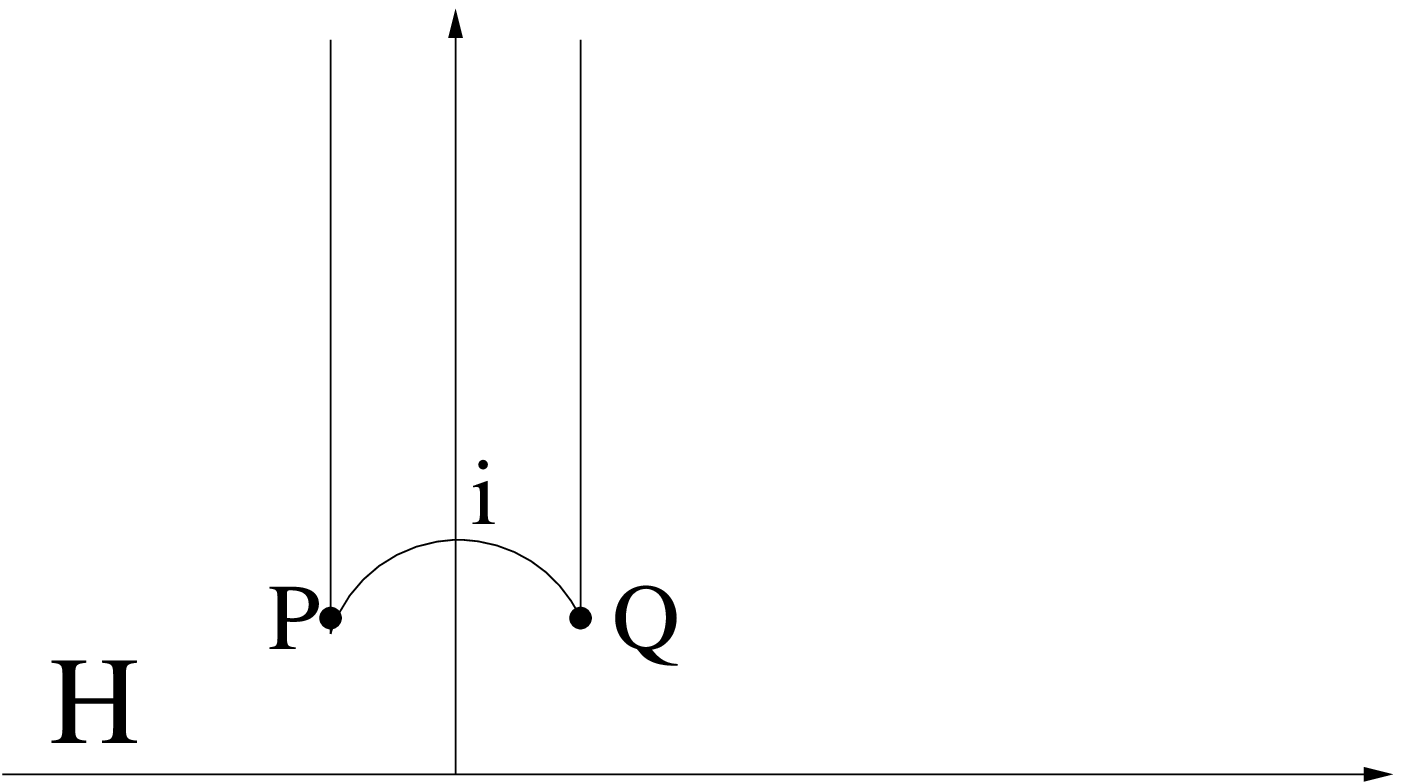}\\
\vspace*{.5mm}
\refstepcounter{diagramm}{\it Figure \arabic{diagramm}}:
{\it Fundamental domain of $\slzwei(\ZZ)$.}
\label{fd1}
\end{center}

The surface $\HH/\slzwei(\ZZ)$ is obtained by identifying the vertical 
edges $P\infty$ and $Q\infty$ via $T$ and the  edge PQ with itself
(with fixed point $i$) via $S$.\\
For an arbitrary subgroup $\Gamma$ of $\slzwei(\ZZ)$ of finite index
we obtain a fundamental
domain as a union of translates of the triangle $\Delta$: for each coset
$\overline{A}$ we take the triangle $A(\Delta)$, where $A$ is a representative
of the coset. The identification of the edges is given by the elements 
in $\Gamma$. Gluing the edges gives the quotient surface $\HH/\Gamma$, 
filling in the 
cusps leads to a closed Riemann surface endowed with a triangulation.
We draw stylized pictures of the fundamental domains that indicate the
triangles (see Figure \ref{fdl23} and \ref{dcurve}). 
The  triangles are labeled with a coset representative,  edges
that are identified are labeled with the same letter and vertices that
are identified with the same number. Vertices that come from cusps 
(i.e. points at $\infty$) are marked with $\bullet$. \\
In particular, one can read off from these stylized pictures the genus and
the number of cusps of the quotient surface $\HH/\Gamma$.

\subsubsection*{Two examples: 
the origami L(2,3) and 
the origami D}

{\bf The origami L(2,3):}\\
In \cite[Example 3.5]{sh1} the Veech group is calculated as follows:

\[\veech(L(2,3)) \; = \;
   <\begin{pmatrix}1&3\\0&1\end{pmatrix},
    \begin{pmatrix}1&0\\2&1\end{pmatrix}, \begin{pmatrix}-1&3\\-2&5\end{pmatrix},
    \begin{pmatrix} 3&-5\\2&-3\end{pmatrix},
    \begin{pmatrix} -1&0\\0&-1\end{pmatrix}>.\]

More precisely, one obtains the generators presented as products
of $S$ and $T$ as well as a list of coset representatives.
\begin{itemize}
\item
List of generators:
\begin{eqnarray*}
\begin{pmatrix}1&3\\0&1\end{pmatrix} \,= \,
T^3, \quad
\begin{pmatrix}-1&3\\-2&5 \end{pmatrix} \,=\,
TST^2ST^{-1}T^{-1},\quad
\begin{pmatrix}1&0\\2&1 \end{pmatrix} \,=\,
TSTST^{-1}S,\\[2mm]
\begin{pmatrix}3&-5\\2&-3 \end{pmatrix} \,=\,
T^2STST^{-1}S^{-1}T^{-2}, \quad
\bpm -1&0\\0&-1 \epm \,=\, -I 
\hspace*{3.6cm}
\end{eqnarray*}
\item 
List of representatives:
\[I,\,\; T,\,\; S,\,\; T^2, \,\;TS,\,\;
 ST,\,\; T^2S,\,\;
 TST,\,\; T^2ST
\]
\end{itemize}

Hence, $\Gamma(L(2,3))$ is a subgroup of index $9$ in $\slzwei(\ZZ)$.

The stylized picture of the quotient $\HH/\Gamma(L(2,3))$ is determined
in \cite[Example 3.6]{sh1} and is shown here in Figure
\ref{fdl23}.
\vspace*{-5mm}
\begin{center}
\includegraphics[scale = .5]{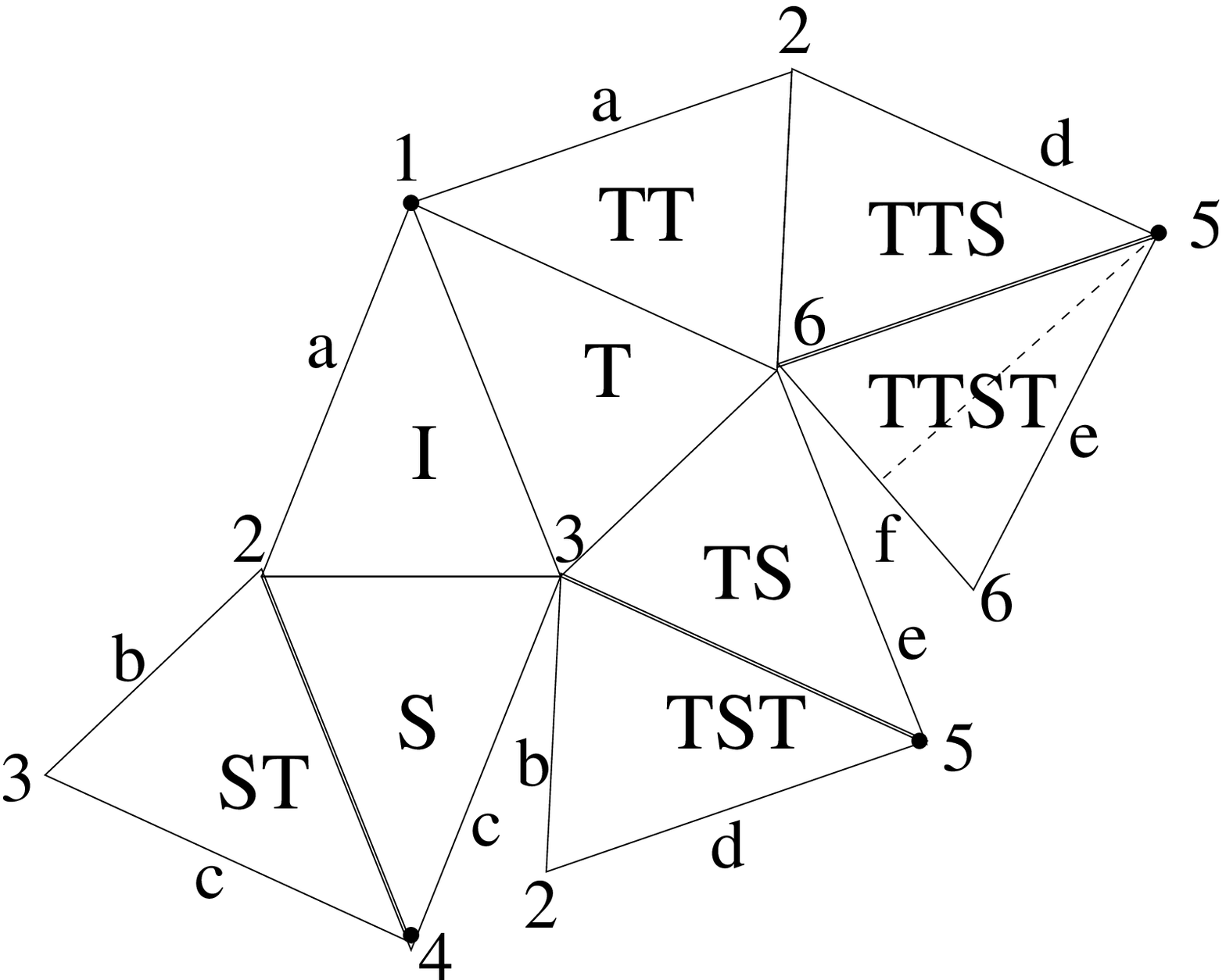}\\
%\end{center}
%\vspace*{-10mm}
%\begin{center}
\refstepcounter{diagramm}{\it Figure \arabic{diagramm}}:
{\it Fundamental domain of $\Gamma(L(2,3))$.}
\label{fdl23}
\end{center}

From this one can read off that the genus of the quotient 
$\HH/\Gamma(L(2,3))$ is $0$ and that it has 3 cusps, namely the vertices
labeled by 1,4 and 5.
It follows in particular that the corresponding Teichm\"uller curve
has genus 0.\\

{\bf The origami D:}\\
The Veech group of the origami $D$ is calculated in
\cite[Section 7.3.2]{gabidiss}. It  has index $24$ in $\slzwei(\ZZ)$
and the following generators:
%\footnote{The data were calculated with a computer program
%that implements the algorithm presented in Section \ref{algo}.}

\[\begin{array}{lcllcl}
A_0 &=&\begin{pmatrix}  -1& 0\\  0&-1\end{pmatrix} \,=\, -I,&
A_1 &=& \begin{pmatrix}  1&3 \\  0& 1\end{pmatrix} \,=\, T^3,\\[5mm]
A_2 &=& \begin{pmatrix}  1&0 \\ -6& 1\end{pmatrix} \,=\,ST^6S^{-1},&
A_3 &=& \begin{pmatrix} -7&16\\ -4& 9\end{pmatrix} \\[5mm]
    &  &                                                    &
    &\phantom{:}=&      (T^2S)T^4(T^2S)^{-1}.\\[5mm]
A_4 &=& \begin{pmatrix} -3&4\\ -4& 5 \end{pmatrix} \,=\,(TS)T^4(TS)^{-1},&
A_5 &=& \begin{pmatrix}  -9&5\\-20&11\end{pmatrix} \\[5mm]
    &  &                                                    &
    &\phantom{:}=&                   (TST^2S)T^5(TST^2S)^{-1},\\[5mm]
A_6 &=& \begin{pmatrix} 7&2\\-18&-5\end{pmatrix} \,=\,(ST^3S)T^2(ST^3S)^{-1},
 \hspace*{-1.5cm}\\[5mm]
\end{array}\]

The following is a system of cosets representatives:
\begin{eqnarray*}
&&I\,,\ T\,,\ S\,,\  T^2\,,\  TS\,,\  ST\,,\  T^2S \,,\ TST\,,\  ST^2\,,\ STS \, ,\,
 T^2ST \,,\  TST^2 \,,\\\
&& ST^5 \,,\ ST^3  \,,\  T^2S \,,\ TST^3 \,,\  TST^2S \,,\  ST^4 \,,\ ST^3S \,,\
   TST^2ST^{-1} \,,\ \\
&& TST^2ST^{-2} \,,\ TST^2ST^{-3} \;;\  TST^2ST^{-4} \,,\  ST^3ST
\end{eqnarray*}

The corresponding origami curve $C(D)$ has genus $0$. It
is shown with its natural triangulation 
in  Figure  \ref{dcurve}. It has six cusps, namely $C_1$, $C_2$, $C_3$, 
$C_4$, $C_5$ and $C_6$. 
 \\

\hspace*{3mm}\includegraphics[scale = .65]{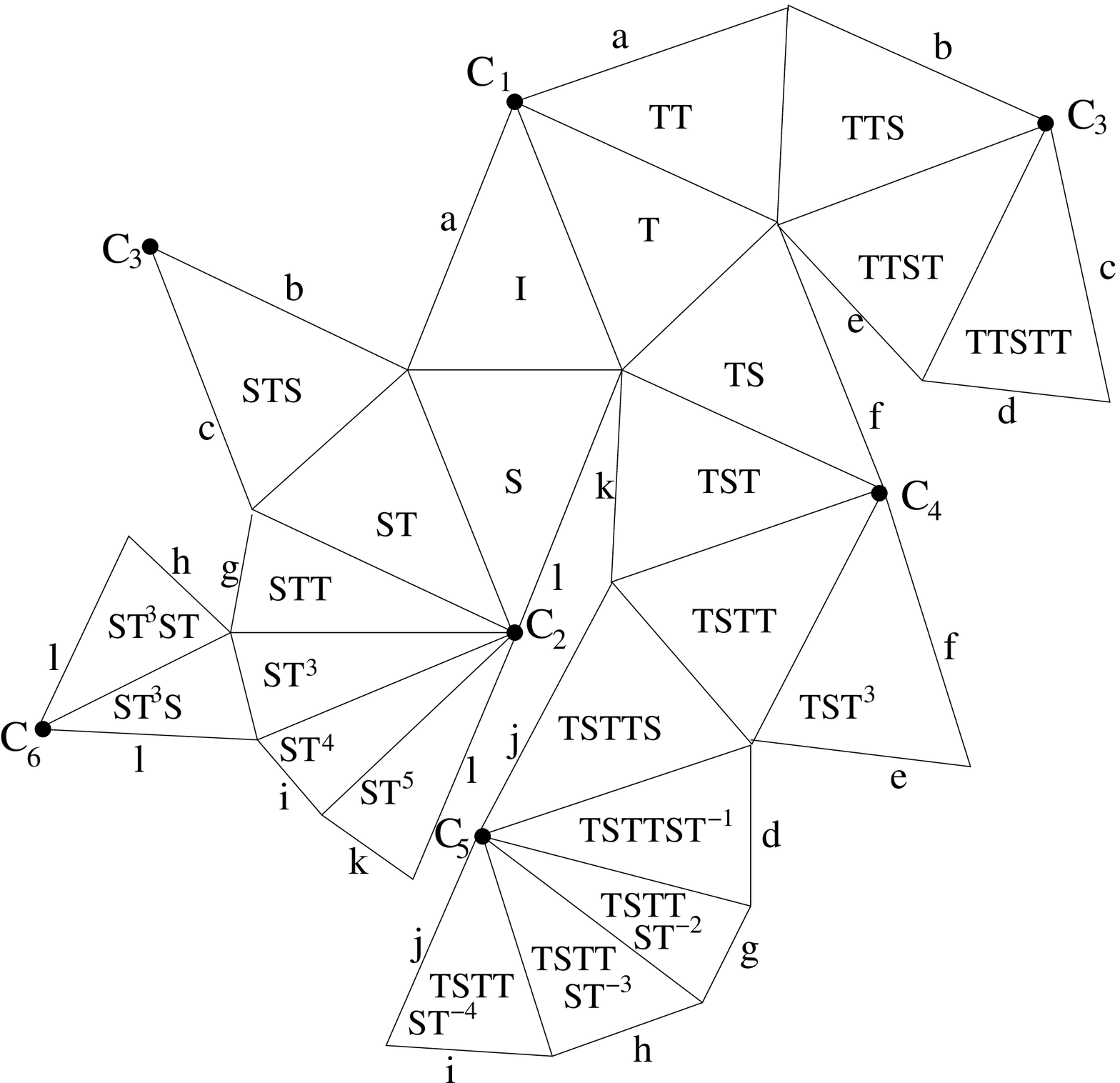}

\begin{center}
\refstepcounter{diagramm}{\it Figure \arabic{diagramm}}:
 {\it \begin{tabular}[t]{l}
   The origami curve to $D$. \\
   \end{tabular}
}
\label{dcurve}
\end{center}

\section{Veech groups that are non congruence groups}\label{noncgsection}

Theorem \ref{cgs} implies that there are many congruence groups which
are Veech groups. How about non congruence groups? In this section
we will see that the Veech groups for the two examples, the origami
$L(2,3)$ and the origami $D$, studied in
the last paragraph are both non
congruence groups. Furthermore, 
we give a construction that produces for both of them an infinite
sequence of origamis whose Veech group is a non congruence group.
We use this in order to prove our main theorem.\\

An other generalization of the example $L(2,3)$ was given by Hubert
and Leli\`evre in \cite{hl2}, where they show for certain ``L-shaped''
origamis or square-tiled surfaces, how they are called there, 
that their Veech groups are non congruence groups. These
surfaces are all of genus 2, hence it follows that there are infinitely
many origamis of genus 2 whose Veech group is a non congruence group. \\

Recall that a group is a congruence group, whose level is a divisor
of $n$,  
if and only if it contains the principal congruence group
\[\Gamma(n) \; = \; \{\bpm a&b\\c&d\epm \equiv \bpm 1&0\\0&1\epm
  \mod n\}
  \; = \;\; \mbox{kernel}(\mbox{proj}:\slzwei(\ZZ) \to \slzwei(\ZZ/n\ZZ))
\]

In \cite[Proposition 3.8]{sh1} it was shown using a proof of Stefan 
K\"uhnlein that the Veech group of $L(2,3)$ is a non congruence group. The
basic tool for this is the general level that is defined for any 
subgroup $\Gamma$ of $\slzwei(\ZZ)$ as follows: For each cusp we define 
its {\em amplitude}
to be the smallest natural number $n$ such that there is an element 
of $\Gamma$ conjugated in $\slzwei(\ZZ)$
to the matrix
\[\bpm 1 &n\\ 0&1 \epm\]
which fixes the cusp. Observe that this is equal to the number of triangles
around the vertex that represents the cusp in our stylized picture
of the quotient surface (see Figures \ref{fdl23} and \ref{dcurve}). 
The {\em general level} of $\Gamma$ 
is the least common multiple of the amplitudes
of all its cusps. 
A theorem of Wohlfahrt \cite[Theorem 2]{wo} states that the level
and the general level of a congruence group coincide. \\
The amplitude of the three cusps of $\HH/\Gamma(L(2,3))$ labeled
with 1, 4 and 5 in Figure \ref{fdl23} is 
3, 2 and 4 respectively. Hence, the general level of $\Gamma(L(2,3))$ 
is 12. Then it is shown in the proof that 
$\Gamma(L(2,3))$ does not contain $\Gamma(12)$ which gives the contradiction.\\

The same method can be used in order to show that $\Gamma(D)$ is a non 
congruence group. We here carry out the proof for it. Observe 
from Figure \ref{dcurve} that the six cusps $C_1$, \ldots,  $C_6$ have 
the amplitude 3, 6, 4, 4, 5 and 2, respectively. Thus the general level
is 60.

\begin{prop}
The Veech group $\Gamma(D)$  is a non congruence group.
\end{prop}

\begin{proof}
Suppose that $\Gamma = \Gamma(D)$ is a congruence group.
Since the general level of $\Gamma$ is $60$, we have  by the 
theorem of Wohlfahrt mentioned above, 
that $\Gamma(60)$ is a subgroup of $\Gamma$.\\

We will use the following facts, which can be checked e.g. 
in Figure \ref{dcurve}:
\[\begin{array}{c}
  A_1 = \begin{pmatrix}  1&3 \\ 0& 1\end{pmatrix}  \in \Gamma, \quad
  A_6 = \begin{pmatrix}  7& 2\\-18&-5\end{pmatrix} 
   \in \Gamma \quad \mbox{ and }
  T = \begin{pmatrix}1&1\\0&1\end{pmatrix} \notin \Gamma
  \end{array}\]
In order to verify this in Figure \ref{dcurve}, use that
\[A_1 = T^3  \mbox{ and } 
  A_6 = S^{-1}T^2S^{-1}T^{-1}S^{-1}TS^{-1}T^{-3}S^{-1}.\]
We will find an element in 
$\Gamma$ whose projection to
$\slzwei(\ZZ/60\ZZ)$ is equal to that of $T$, which gives
us the desired contradiction.\\
Recall that
\[  \slzwei(\ZZ/60\ZZ) \cong \slzwei(\ZZ/4\ZZ) \times \slzwei(\ZZ/3\ZZ)
                                            \times \slzwei(\ZZ/5\ZZ). \]
We identify in the following these two groups.
Furthermore we denote by $p_4$, $p_3$, $p_5$ and $p_{60}$
the projection from $\slzwei(\ZZ)$  to $\slzwei(\ZZ/4\ZZ)$,
$\slzwei(\ZZ/3\ZZ)$, $\slzwei(\ZZ/5\ZZ)$ and $\slzwei(\ZZ/60\ZZ)$, respectively.
Then $p_{60} = p_4 \times p_3\times p_5$.\\
We have
\begin{eqnarray*}
p_{60}(A_1) &=& (\begin{pmatrix} 1&3 \\ 0& 1\end{pmatrix},
            \begin{pmatrix}  1&0 \\ 0& 1\end{pmatrix},
            \begin{pmatrix}  1&3 \\ 0& 1\end{pmatrix}) \;\;\quad \mbox{ and}\\
p_{60}(A_6) &=& (\begin{pmatrix}  3& 2\\2&3\end{pmatrix},
            \begin{pmatrix}  1& 2\\0&1\end{pmatrix},
            \begin{pmatrix}  2& 2\\2&0\end{pmatrix})
\end{eqnarray*}
The order of $p_4(A_1)$ in $\slzwei(\ZZ/4\ZZ)$ is $4$,
the order of $p_3(A_1)$ in $\slzwei(\ZZ/3\ZZ)$ is $1$ and
the order of $p_5(A_1)$ in $\slzwei(\ZZ/5\ZZ)$ is $5$. We also say:
{\it The order of $p_{60}(A_1)$ is $(4,1,5)$.}
Since \; $7 \equiv 3 \mod 4$ \;  and \; $7 \equiv 2 \mod 5$ \; we have
\begin{equation}\label{Aeins}
p_{60}(A_1^7) =
           (\begin{pmatrix} 1&1 \\ 0& 1\end{pmatrix},
            \begin{pmatrix}  1&0 \\ 0& 1\end{pmatrix},
            \begin{pmatrix}  1&1 \\ 0& 1\end{pmatrix})
\end{equation}
Furthermore:
\[p_{60}(A_6^2) =(\begin{pmatrix} 1&0 \\ 0& 1\end{pmatrix},
            \begin{pmatrix}  1&1 \\ 0& 1\end{pmatrix},
            \begin{pmatrix}  3&4 \\ 4& 4\end{pmatrix}) \]
and with the same notation as above {\it $p_{60}(A_6^2)$ has the order
$(1,3,5)$}. Thus
\begin{equation}\label{Asieben}
p_{60}(A_6^{20}) =(\begin{pmatrix} 1&0 \\ 0& 1\end{pmatrix},
                 \begin{pmatrix}  1&1 \\ 0& 1\end{pmatrix},
                 \begin{pmatrix}  1&0 \\ 0& 1\end{pmatrix}).
\end{equation}
From (\ref{Aeins}) and (\ref{Asieben}) it follows that
\[p_{60}(A_6^{20}\cdot A_1^7) \;=\;
  (\begin{pmatrix} 1&1 \\ 0& 1\end{pmatrix},
                 \begin{pmatrix}  1&1 \\ 0& 1\end{pmatrix},
                 \begin{pmatrix}  1&1 \\ 0& 1\end{pmatrix})
   \;=\; p_{60}( \begin{pmatrix}  1&1 \\ 0& 1\end{pmatrix})
   \;=\; p_{60}(T) \]
But $A_6^{20}\cdot A_1^7 \in \Gamma$ and $T \notin \Gamma$,
thus $\Gamma(60) = \ker(p_{60})$ cannot be contained in $\Gamma$.
Therefore, $\Gamma$ cannot be a congruence group of level $60$. Contradiction!
\end{proof}

\subsubsection*{Sequences of origamis with non congruence Veech groups}

Starting from the origamis $L(2,3)$ and $D$ we will define respectively
a sequence $O_n$, such that for each $n \in \NN$ the Veech group $\Gamma(O_n)$
again is a non congruence group. 
The basic idea is to ``copy and paste'': we will cut the
origami along a segment, take $n$ copies of it and glue them along the cuts.\\

In Figure \ref{locn} we show the origami $O_n$ for $L(2,3)$:

\begin{center}

\setlength{\unitlength}{1cm}
\hspace*{1cm}
\begin{picture}(13,2.3)
% Die Kaestchen des Origamis:
%\put(-2.7,1.2){ $\mbox{Loc}(n) =$}
\put(0,0){\framebox(1,1){1}}
\put(1,0){\framebox(1,1){3}}
\put(2,0){\framebox(1,1){4}}
\put(0,1){\framebox(1,1){2}}

\put(3,0){\framebox(1,1){5}}
\put(4,0){\framebox(1,1){7}}
\put(5,0){\framebox(1,1){8}}
\put(3,1){\framebox(1,1){6}}

\put(6.3,0.5){\ldots}

\put(7.3,0){\framebox(1,1){4n-7}}
\put(8.3,0){\framebox(1,1){4n-5}}
\put(9.3,0){\framebox(1,1){4n-4}}
\put(7.3,1){\framebox(1,1){4n-6}}
\put(10.3,0){\framebox(1,1){4n-3}}
\put(11.3,0){\framebox(1,1){4n-1}}
\put(12.3,0){\framebox(1,1){4n}}
\put(10.3,1){\framebox(1,1){4n-2}}

% labels of the vertices:
\put(-.1,-.1){\Large $\bullet$}
%\put(-.2,-.4){$P_1$}
\put(-.1, .88){\Large $\bullet$}
\put(-.1,1.88){\Large $\bullet$}
\put( .88,-.1){\Large $\bullet$}
\put( .88, .88){\Large $\bullet$}
\put( .88,1.88){\Large $\bullet$}
\put(1.88,-.1){\Large\bf $\circ$}
\put(1.88, .88){\Large $\circ$}

\put(2.88,-.1){\Large $\bullet$}
\put(2.88, .88){\Large $\bullet$}
\put(2.88,1.88){\Large $\bullet$}
\put(3.88,-.1){\Large $\bullet$}
\put(3.88, .88){\Large $\bullet$}
\put(3.88,1.88){\Large $\bullet$}
\put(4.86,-.1){\Large\bf $\circ$}
\put(4.86, .88){\Large $\circ$}
\put(7.18,-.1){\Large $\bullet$}
\put(7.18, .88){\Large $\bullet$}
\put(7.18,1.88){\Large $\bullet$}
\put(8.18,-.1){\Large $\bullet$}
\put(8.18, .88){\Large $\bullet$}
\put(8.18,1.88){\Large $\bullet$}
\put(9.17,-.1){\Large\bf $\circ$}
\put(9.17, .86){\Large $\circ$}

\put(10.18,-.1){\Large $\bullet$}
\put(10.18, .88){\Large $\bullet$}
\put(10.18,1.88){\Large $\bullet$}
\put(11.18,-.1){\Large $\bullet$}
\put(11.18, .88){\Large $\bullet$}
\put(11.18,1.88){\Large $\bullet$}
\put(12.17,-.1){\Large\bf $\circ$}
\put(12.17, .86){\Large $\circ$}

\put(13.18, -.1){\Large $\bullet$}
\put(13.18, .88){\Large $\bullet$}

\end{picture}

\end{center}
\begin{center}
\refstepcounter{diagramm}{\it Figure \arabic{diagramm}}:
{\it $n$ copies of $L(2,3)$. Opposite edges are glued.}
\label{locn}\\
\end{center}

Using the description of an origami by a pair of permutations from Section 
\ref{origamis}, $O_n$ is given as:
\[\sigma_a = (1\;3\;4\;5\;7\;8\;9\;11\;12\;\ldots\; 4n-3\; 4n-1\; 4n),\quad      \sigma_b = (1\;2)(5\;6)\ldots(4n-3\;4n-2).\] 
Observe that the genus of $O_n$ is $n+1$ 
and it has $2n$ cusps: 
$n$ of order $3$ (all $n$ marked by $\bullet$ in Figure \ref{locn}), 
and $n$ of order $1$ (all $n$ marked by $\circ$ in  Figure \ref{locn}). \\

Finally, we want to present the origami $O_n$ by the finite index subgroup 
$H_n = \pi_1(\Xstern)$ of $F_2$, 
that corresponds to $O_n$ by Remark \ref{corespond}.\\ 
Recall from
Example \ref{fgex} that for $O_1 = L(2,3)$, we obtain
the free group of rank 5: \\
\[U \,=\, H_1 \,=\, <g_1 = x^3,\; g_2 = xyx^{-1},\; 
      g_3 = x^2yx^{-2},\; g_4 = yxy^{-1},\; 
      g_5 = y^2> \,= \,F_5.\]
The group $H_n$ is obtained  as 
as follows:
\[H_n \;=\; <g_1^{\,n},\;\; g_1^i\,g_j\,g_1^{-i} \; 
   \in F_5\,|\;\; i \in \{0, \ldots, n-1\}
 \mbox{ and } j \in \{2,\ldots, 5\}>\]

In Figure \ref{dn}, we show the origami $D_n$:

\begin{center}
\setlength{\unitlength}{1cm}
\hspace*{1cm}
\begin{picture}(12,3.2)
% Die labels of the vertices des Origamis:
%\put(-2.7,1.2){ $\mbox{Loc}(n) =$}
\put(0,0){\framebox(1,1){1}}
\put(1,0){\framebox(1,1){2}}
\put(2,0){\framebox(1,1){3}}
\put(0,1){\framebox(1,1){4}}
\put(0,2){\framebox(1,1){5}}

\put(3,0){\framebox(1,1){6}}
\put(4,0){\framebox(1,1){7}}
\put(5,0){\framebox(1,1){8}}
\put(3,1){\framebox(1,1){9}}
\put(3,2){\framebox(1,1){10}}
 
\put(6.3,0){\ldots\ldots\ldots}

\put(8,0){\framebox(1,1){5n-4}}
\put(9,0){\framebox(1,1){5n-3}}
\put(10,0){\framebox(1,1){5n-2}}
\put(8,1){\framebox(1,1){5n-1}}
\put(8,2){\framebox(1,1){5n}}

% labels:
\put(1.4,-.4){$a_1$}
\put(2.4,1.2){$a_1$}
\put(1.4,1.2){$b_1$}
\put(2.4,-.4){$b_1$}

\put(4.4,-.4){$a_2$}
\put(5.4,1.2){$a_2$}
\put(4.4,1.2){$b_2$}
\put(5.4,-.4){$b_2$}

\put(9.4,-.4){$a_{n}$}
\put(10.4,1.2){$a_{n}$}
\put(9.4,1.2){$b_{n}$}
\put(10.4,-.4){$b_{n}$}

% labels Punkte:
\put(-.1,-.1){\Large $\bullet$}
\put( .9,-.1){\Large $\bullet$}
\put(2.9,-.1){\Large $\bullet$}
\put(1.9, .9){\Large $\bullet$}
\put(-.1,2.9){\Large $\bullet$}
\put( .9,2.9){\Large $\bullet$}
\put(-.12, 1.88){\Large $\circ$}
\put( .88, 1.88){\Large $\circ$}
\put(-.17,  .69){\LARGE \bf *}
\put( .83,  .69){\LARGE \bf *}
\put(1.83,  -.31){\LARGE \bf *}
%\put(2.83,  .69){\LARGE \bf *}

\put(2.9,-.1){\Large $\bullet$}
\put(3.9,-.1){\Large $\bullet$}
\put(5.9,-.1){\Large $\bullet$}
\put(4.9, .9){\Large $\bullet$}
\put(2.9,2.9){\Large $\bullet$}
\put(3.9,2.9){\Large $\bullet$}
\put(2.88, 1.88){\Large $\circ$}
\put(3.88, 1.88){\Large $\circ$}
\put(2.83,  .69){\LARGE \bf *}
\put(3.83,  .69){\LARGE \bf *}
\put(4.83,  -.31){\LARGE \bf *}
\put(5.83,  .69){\LARGE \bf *}

\put(7.9,-.1){\Large $\bullet$}
\put(8.9,-.1){\Large $\bullet$}
\put(10.9,-.1){\Large $\bullet$}
\put(9.9, .9){\Large $\bullet$}
\put(7.9,2.9){\Large $\bullet$}
\put(8.9,2.9){\Large $\bullet$}
\put(7.88, 1.88){\Large $\circ$}
\put(8.88, 1.88){\Large $\circ$}
\put(7.83,  .69){\LARGE \bf *}
\put(8.83,  .69){\LARGE \bf *}
\put(9.83,  -.31){\LARGE \bf *}
\put(10.83,  .69){\LARGE \bf *}

\end{picture}

\end{center}
\begin{center}
%\refstepcounter{diagramm}{\it Figure \arabic{diagramm}}:
%{\it The origami $L_n$:
\refstepcounter{diagramm}{\it Figure \arabic{diagramm}}:
{\it 
$n$ copies of $D$. 
\begin{tabular}[t]{l}
Edges with the same label or\\ 
unlabeled opposite edges are glued.
\end{tabular}}
\label{dn}\\
\end{center}
 The pair of
permutations describing $D_n$ is: 
\begin{eqnarray*}
  \sigma_a &=& (1 \; 2\;3 \;\, \; 6 \; 7\; 8 \;\, \ldots \;\,
              5n-4 \; 5n-3 \; 5n-2), \\
  \sigma_b &=& (1\;4\;5)(6\;9\;10)\ldots
     (5n-4\;5n-1\;5n)(2\;3)(7\;8)\ldots(5n-3\;5n-2)
\end{eqnarray*}
The genus of $D_n$ is $2n$ and it has $n+2$ cusps: 
2 of order $2n$ (marked as $\bullet$ and $\star$)
and $n$ of order $1$ (all $n$ marked by $\circ$).\\

Again, we present $O_n$ by the corresponding finite
index subgroup $H_n$ of $F_2$. We have from Example \ref{fgex}
that $U = H_1 = F_6$, the free group of rank 6:
\[
U \, = \, 
  <g_1' =  x^3,\; g_2' = xyx^{-2},\; 
   g_3' = x^2yx^{-1},\; g_4' = yxy^{-1},\; 
   g_5' = y^2xy^{-2},\; g_6' = y^3> \\
      = F_6
\]
And similarly as above, we obtain:
\[H_n \;=\; <g_1'^{\,n},\;\; g_1'^i\,g_j'\,g_1'^{-i} \; 
   \in F_6\,|\;\; i \in \{0, \ldots, n-1\}
 \mbox{ and } j \in \{2,\ldots, 6\}>\]

We will see in the following that for both sequences all Veech groups
$\Gamma(O_n)$ are non congruence groups. More precisely, we will
show:

\begin{prop}\label{hp}
For both sequences $O_n$ the following holds:
\begin{itemize}
\item 
$\Gamma(O_n) \subseteq \Gamma(O_1)$, 
which is for both sequences a non
congruence group. 
\item More generally one has:\\
$n$ divides $m$ \;$\Rightarrow$ \;
$\Gamma(O_m) \subseteq \Gamma(O_n)$.
\item
Different origamis in one sequence have different Veech groups, 
i.e.:\\
$\Gamma(O_n) \neq \Gamma(O_m)$ for $n \neq m$.
\end{itemize}
\end{prop}

To prove this, let us 
detect that we are in the following more general
setting.\\

{\bf Setting A:}
\begin{itemize}
\item
Let $U$ be a finite index subgroup of $F_2$. Then
$U$ is a free group of rank $k$ for some $k \geq 2$, i.e.
\[U \; = \; <g_1, \ldots, g_k> \; = \; F_k\]
\item
Let $\alpha: F_k \to \ZZ$ be the projection 
$w \mapsto \sharp_{g_1} w$\\
where $\sharp_{g_1} w$ is the number of $g_1$ in the word 
$w = w(g_1,\ldots, g_k)$ with $g_1^{-1}$ counted as $-1$.
\item 
Let $H_n$ be the kernel of $p_n \circ \alpha$,
where $p_n:\ZZ \to \ZZ/n\ZZ$ is the natural projection, i.e.
\[H_n \;=\; <g_1^{\,n},\;\; g_1^i\,g_j\,g_1^{-i} \; \in F_k\,|
  \;\; i \in \{0, \ldots, n-1\}
 \mbox{ and } j \in \{2,\ldots, k\}>.\]
\item
Finally, let $H_{0}$ be the kernel of $\alpha$, i.e.:
\[H_{0} \;=\; \bigcap_{n\in\NN} H_n \;=\; <<g_2, \ldots, g_k>>_U,\]
is the normal subgroup in $U$ generated by $g_2$, \ldots, $g_k$. 
\end{itemize}

Observe that we are in this setting with\\[1mm]
$U \,=\, \pi_1(\Xstern) 
   \,=\, <x^3,\; xyx^{-1},\; x^2yx^{-2},\; 
  yxy^{-1},\; y^2>$
\; for the origami $L(2,3)$ and\\[1mm]
$U \,=\, \pi_1(\Xstern) 
   \,=\, <x^3,\; xyx^{-2},\; x^2yx^{-1},\; yxy^{-1},\; y^2xy^{-2},\; y^3>$ 
\; for the origami $D$.\\

In order to prove the properties in Proposition \ref{hp},
we will need 
that $U$ fulfills the following a bit technical condition:\\[2mm]
{\bf Property B:} 
Let $U \,=\, <g_1, \ldots, g_k>$ \, $(k \geq 2)$
be as above a finite index subgroup of $F_2$
of rank $k$ and $\{w_i\}_{i \in I}$  a system
of coset representatives with $w_1 = \id$. 
Suppose that $U$ has the following property:
\begin{equation*}\label{blubs}
\forall \, j \in I - \{1\}:\;\; w_j\,<<g_2,\ldots, g_k>>_U\,w_j^{-1} \not\subseteq U.
\end{equation*}

One can check by hand that for both origamis, $L(2,3)$ and $D$, this
property is fulfilled. In this setting we obtain the following
conclusions.

\begin{prop} \label{cr}
Let  $n \in  \NN \, \cup \, \{0\}$. Let $U$ be
a finite index subgroup of $F_2$ fulfilling
property B. With the
notations from Setting A, we have:
\begin{enumerate}
\item[a)]
The normalizer of $H_n$ in $F_2$ is equal to $U$:\; 
{\em Norm}$_{F_2}(H_n) \,=\, U$
\item[b)]
$\emstabaut(H_n) \; \subseteq \; \emstabaut(U)
 \; \stackrel{\mbox{\scs \em Def}}{=}\; \GGG$ 
\item[c)]
Recall that $U = F_k$, the free group in $k$ generators.\\
Let $\beta_n: \emaut(F_k) \to \emglk(\ZZ/n\ZZ)$ be
the natural projection.\\[1mm]
Then
$\emstabaut(H_n)$ is equal to
\[\beta_n^{-1}(\{A = (a_{i,j})_{1 \leq i,j \leq k} \in \, \emglk(\zn)|
         \; a_{1,2} = \ldots = a_{1,k} = 0\}) \;\cap \; \GGG.\]
Here we use the notation $\ZZ/(0\ZZ) = \ZZ$ thus $\beta_0$ is
the natural projection $\emaut(F_k) \to \mbox{\em GL}_k(\ZZ)$.
\end{enumerate}
\end{prop}

\begin{proof}\hspace*{1cm}\\
{\bf a)}\\
By definition $H_n$ is normal in $U$, i.e.\ 
$U \subseteq \norm_{F_2}(H_n)$.\\ 
Let now $w$ be an
element of $F_2 \backslash U$. Hence, $w = w_j\cdot u$ for some 
$j \in I-\{1\}$, $u \in U$. 
By Property B, there exists some $h_{0} \in \; 
<<g_2,\ldots, g_k>>_U \;=\,H_{0}$, such that
$w_jh_{0}w_j^{-1} \not\in U$. Therefore we have
$w(u^{-1}h_{0}u)w^{-1} \not\in U$. But 
$u^{-1}h_{0}u \in H_{0} \subseteq H_{n}$, since $H_{0}$ 
is normal in $U$. This shows that $w \not\in \norm_{F_2}(H_n)$.\\[2mm]
{\bf b)}\\
This follows from a), since for a subgroup $H$ of $F_2$ in general holds:\\
$\stabaut(H) \subseteq \stabaut(\norm_{F_2}(H))$, see e.g. 
\cite[Remark 3.1]{sh2}.\\[2mm]
{\bf c)}\\
Define $M =  \{A = (a_{i,j})_{1 \leq i,j \leq k} \in \, \glk(\zn)|
         \; a_{1,2} = \ldots = a_{1,k} = 0\}$.\\
Let $\gamma \in \GGG$. We have to show that $\gamma(H_n) = H_n$
if and only if $\beta_n(\gamma) \in M$.\\
Let furthermore $p_n^k: F_k \to (\ZZ/n\ZZ)^k$ be the natural projection.\\
Consider the following commutative diagram:
\[ \xymatrix{
   H_n \ar[d]^{p_n^k} & \subseteq &  
   F_k  \ar[r]^{\gamma} \ar[d]^{p_n^k}& 
   F_k \ar[d]^{p_n^k}& 
   \supseteq & H_n \ar[d]^{p_n^k}\\
   \Hquer_n = p_n^k(H_n) &  \subseteq & (\ZZ/n\ZZ)^k \ar[r]^{\beta_n(\gamma)} &
   (\ZZ/n\ZZ)^k & \supseteq  & p_n^k(H_n) = \Hquer_n
   }
\]
Since $p_n^k$ is surjective and $H_n$ is the full preimage of 
$\Hquer_n = p_n^k(H_n)$,
it follows that 
$\gamma(H_n) = H_n$ if and only if 
$\beta_n(\gamma)(\Hquer_n) = \Hquer_n$.\\[2mm]
Observe finally that:
\begin{eqnarray*}
\Hquer_n &=& \{(0,x_2, \ldots, x_k) \in (\ZZ/n\ZZ)^k\} \quad \mbox{ and }\\[2mm]
\stab_{\,\mbox{\scs GL}_k(\ZZ/n\ZZ)}(\Hquer_n) 
    &=& \{\begin{array}[t]{l}
         A = (a_{i,j})_{1 \leq i,j \leq k} \in \glk(\ZZ/n\ZZ) | \; \,\\
        (y_1, \ldots, y_k) = A\cdot(0,x_2, \ldots, x_k) \;
        \hspace*{\fill}\Rightarrow\; y_1 = 0 \;\}
        \end{array}\\
    &=& \{A = (a_{i,j}) \in \glk(\ZZ/n\ZZ)|
           \; a_{1,2} \,=\, \ldots \,= \,a_{1,k} \,=\, 0\} \\
     &=&  M.
\end{eqnarray*}
\end{proof}

Theorem \ref{char} suggests the following notation.

\begin{defn} Let $U$ be a subgroup of $F_2$.\\ 
With 
$\hat{\beta}: \emautplus(F_2) \to \emslzwei(\ZZ)$ as in Theorem \ref{char}, we
define
\[\Gamma(U) = \hat{\beta}(\emstabaut(U))\]
and call  $\Gamma(U)$ the {\em Veech group} of $U$.
\end{defn}

We now obtain from Proposition \ref{cr} the following conclusions.

\begin{cor} \label{qay} \label{Hinftyschnitt}
Suppose that we are in the same situation as in Proposition \ref{cr}, in 
particular that $U$ is a finite index subgroup of $F_2$
fulfilling property B.  Then we have 
for all $n \in \NN$:
\begin{enumerate}
\item[a)] $\emstabaut(H_{0}) \subseteq \emstabaut(H_n)$ \quad \mbox{ and } 
\quad
      $\Gamma(H_{0}) \subseteq \Gamma(H_{n})$.
\item[b)] If $m \in \NN$ with $n|m$, then:\\[2mm] 
      $\emstabaut(H_{m}) \subseteq \emstabaut(H_n)$ \quad \mbox{ and } \quad
      $\Gamma(H_{m}) \subseteq \Gamma(H_{n})$.
\item[c)] \[\emstabaut(H_{0}) = \bigcap_{n \in \NN} \emstabaut(H_n)
         \;\; \mbox{ and } \quad
        \Gamma(H_{0}) =  \bigcap_{n \in \NN}\Gamma(H_n)  \]
\end{enumerate}
\end{cor}

\begin{proof}\hspace*{1cm}\\
{\bf a) and b):}\\
Let $\gamma \in \GGG$.
By Proposition \ref{cr} we have that
\begin{eqnarray*}
 \forall n \in \NN:\; \gamma \in \stabaut(H_n) &\;\ifff\;& 
   \begin{array}[t]{l}
   \beta_{n}(\gamma) = A = (a_{i,j}) \\[1mm]
    \;\; \mbox{ with } a_{1,2} \,\equiv\, \ldots \,\equiv\, a_{1,k} \,
    \equiv\, 0 \mod n
   \end{array}\\   
  \mbox{and } \quad
   \gamma \in \stabaut(H_{0}) &\;\ifff\;& 
   \begin{array}[t]{l}
   \beta_{0}(\gamma) = A = (a_{i,j}) \\[1mm]
   \;\; \mbox{ with } a_{1,2} \,=\, \ldots \,=\,  a_{1,k} \,=\, 0.
   \end{array}
\end{eqnarray*}
Thus we have for all $n \in \NN$ and for all $m \in \NN$ with $n|m$, that 
\[\stabaut(H_{0}) \subseteq \stabaut(H_{m}) \subseteq \stabaut(H_n).\]
We have in particular by the definition of the Veech group of a subgroup of $F_2$:
\[\Gamma(H_{0}) \;\subseteq\; \Gamma(H_{m}) \;\subseteq\; \Gamma(H_n).\]

{\bf c):}\\
$\subseteq$\;
follows from a). \;
$\supseteq$\;
follows from Remark \cite[Remark 3.1]{sh2}.\\
\end{proof}

We now return to the language of origamis:
Let $O$ be an origami, $U$ the corresponding subgroup of $F_2$. Define for $U$
the subgroups $H_n$ \, ($n\in\NN$) as in Setting A and let 
$O_n$ be the origamis corresponding to the groups $H_n$. \\
By Corollary \ref{qay} and Theorem \ref{char} we obtain immediately
the following result.

\begin{prop}\label{propmulti}
If $U$ has the Property B, then
\[\forall n \in \NN:\; \Gamma(O_n) \subseteq \Gamma(O) \quad\mbox{ and } \quad
   \forall n,m \in \NN: \; 
    n | m \;\Rightarrow\; \Gamma(O_m) \subseteq \Gamma(O_n).\]
In particular, if \, $\Gamma(O)$ is a non congruence group, each $\Gamma(O_n)$
is  a non congruence group. Thus in this case,  we 
obtain infinitely many origamis whose
Veech group is a non congruence group.
\end{prop}

In order to conclude Proposition \ref{hp}, it is now just left to prove 
the last
item. But this follows , since  we have (see \cite{gabidiss}) for both sequences $O_n$, the one
coming from the origami $L(2,3)$ and the one coming from the origami $D$,
that 
\begin{equation}\label{mmm}
\bpm 1&s\\0&1 \epm \in \Gamma(O_n) \;\;\Leftrightarrow\;\; 
   3n \mbox{ divides } s.
\end{equation}

This finishes the proof of Proposition \ref{hp}.\\ 
Furthermore, 
Theorem \ref{main} follows  from
Proposition \ref{hp}. \\

{\bf Remark:}
From Corollary \ref{qay} and (\ref{mmm}) 
it follows that $\Gamma(H_{0})$
has infinite index in $\slzwei(\ZZ)$. Furthermore it is non
trivial, since it contains
\[B_2 = \bpm 1 & 0\\ 2&1 \epm \mbox{ for } L(2,3) \;\;
   \mbox{ respectively } \;\;
  B_3 = \bpm 1 & 0\\ 3&1 \epm \mbox{ for } D.
\]

\end{document}